\documentclass[12pt]{amsart}
\usepackage{latexsym}
\usepackage{psfig}
\usepackage{amsmath,amssymb}
\newtheorem{defn}{Definition}[section]
\newtheorem{lem}[defn]{Lemma}

\newtheorem{prop}[defn]{Proposition}

\newtheorem{cor}[defn]{Corollary}

\newtheorem{them}[defn]{Theorem}

\newenvironment{note}{\noindent{\bf Remark:}}{\\ \par }

\def\co{\colon\thinspace}

\newcommand{\pole}{}

\newcommand{\M}{\mathcal{M}}
\newcommand{\pg}{\partial g}

\author{Anna Klebanov}
\title{Gauss diagrams of ~3-manifolds}
\email{klebanov@techunix.technion.ac.il}
\address{Math Department, Technion, Haifa, 32000, Israel  }

\begin{document}

\begin{abstract}
The paper presents a simple combinatorial method to encode
3-dimensional manifolds, basing on their Heegaard diagrams. The
notion of a Gauss diagram of a ~3-manifold is introduced. We check
the conditions for a Gauss diagram to represent a closed manifold
and a manifold with boundary.
\end{abstract}

\thanks{Partially supported
 by the Israeli Science Foundation, grant 86/01.  }

\maketitle

\section{Introduction}

Orientable ~3-manifolds can be easily described by their Heegaard
splittings. The manifold is presented as a union of two
handlebodies, glued along the boundary by a homeomorphism of a
splitting surface to itself. Any ~3-manifold possesses a Heegaard
splitting, though not a unique one. However, one can pass from one
Heegaard splitting to another by a sequence of some standard
simple operations. See, e.g., \cite{fm}. A standard simple way to
describe a Heegaard splitting of a ~3-manifold is by its  Heegaard
diagram. A Heegaard diagram consists of two families of disjoint
simple closed curves on a splitting surface. The families
intersect, producing an immersed family of curves on the surface.

There is a simple combinatorial way to encode immersed curves in a
plane, namely a Gauss diagram. The curve is presented as a circle
in the plain with double points joined with a chord. One can also
encode the information on the overpasses and underpasses, changing
these diagrams into diagrams for knots and links. This
presentation of plane curves and knot diagrams proved useful for
computational purposes, especially for finite type knot and link
invariants. The same construction can be used to describe immersed
curves on a surface.

A manifold is described by curves on the surface and the curves
are described by Gauss diagrams. It seems natural to introduce a
new notion of Gauss diagrams of ~3-manifolds, which appears to be
a simple combinatorial way to describe ~3-manifolds. In Section
\ref{gd} we introduce the notion of a Gauss diagram corresponding
to a Heegaard diagram of a closed orientable manifold. We show
that the splitting surface and the entire manifold can be
reconstructed from it. The equivalence of Gauss diagrams is
defined and we prove that a closed orientable ~3-manifold is
determined by the equivalence class of its Gauss diagram. Next we
show how to compute the fundamental group and the first homology
group of a manifold directly from a Gauss diagram, and how to
distinguish homology spheres.

A picture looking like a Gauss diagram does not necessary come
from a closed orientable manifold. In Section \ref{abstgd} we
consider all diagrams looking like a Gauss diagram, calling them
abstract Gauss diagrams. We show that this construction allows to
describe oriented manifolds with boundary. Later we formulate the
conditions for a Gauss diagram to represent a closed orientable
manifold and a complement of a knot in some closed orientable
manifold.

The author would like to express her deep gratitude
 to her scientific advisors Prof. Y. Moriah and especially Prof. M. Polyak
 from the Math. Department
 at the Technion, Haifa, Israel for their help, support and patience
 during the research.
 The author also thanks the Technion for financial support.

\section{From ~3-manifolds to Gauss diagrams and back again}
\label{gd} \subsection{Heegaard and Gauss diagrams}

 Let $M$ be a
closed orientable ~3-manifold. A Heegaard splitting
$(H^\pm,\Sigma)$ of $M$ is a union of two handlebodies $H^\pm$,
glued along their boundary $\Sigma$ by a homeomorphism.
 The genus $g$ of $\Sigma$
is called the genus of the Heegaard splitting. One of $H^\pm$ can
be always taken to be standard. The main property of a handlebody
$H$ is that there exists a set of discs $\Delta \subset H$, called
a meridional disc system of $H$, such that $H-\Delta$ is a union
of 3-balls. One can always find a minimal set of this sort, which
means that $H-\Delta$ is a single ball. Since there is a unique
way to fill in the ball, the handlebody can be easily
reconstructed, if a meridional disc system is specified. This
means that it is sufficient to specify the boundaries of the
meridional discs on the surface. Let $\Delta^\pm$ be minimal
meridional disc systems of $H^\pm$ respectively. The boundaries of
these discs are two families $M^\pm$ of disjoint simple closed
curves on $\Sigma$, from which the Heegaard splitting can be
reconstructed by gluing 2-handles along the curves of $M^\pm$ and
filling in 3-balls. The pair $(\Sigma,M^\pm)$ is called a {\em
Heegaard diagram } of $M$.

We think of both families as an immersed family of curves on the
surface. In order to describe them in a more simple, combinatorial
way, one can use Gauss diagrams, in a way similar to the plane
curves. A {\it Gauss diagram } is the immersing collection of
circles with the preimages of each double point connected with a
chord. Each chord is given a sign depending on whether the frame
of tangents to the branches of the curve in this point coincides
with the orientation of the plane. One can also associate each
circle a Gauss word by traveling once around the circle and
recording the chords with signs.

Gauss words were introduced by Gauss in \cite{g}, and the
planarity problem was extensively studied. See, e.g. \cite{lm},
\cite{rr}, \cite{frm}. J.S. Carter proved in \cite{c} that if two
sets of immersed curves fill the surface in the sense that the
complement regions are discs, then they are stable geotopic if and
only if they have equivalent Gauss paragraphs. The notion of Gauss
diagrams for knots and links was introduced by Polyak and Viro in
\cite{pv}. They used Gauss diagrams enhanced with the information
on the overpasses and underpasses. The equivalence classes of
abstract Gauss diagrams were called virtual knots by Kauffman in
\cite{k}\pole,
 since only part
of the diagrams were realizable.

\subsection{From Heegaard diagram to Gauss diagram}
\label{GD}

\begin{defn}
A {\em Gauss diagram} $G=(\M^+,\M^-,h,\varepsilon)$ consists of
two families $\M^\pm$ of $g$ disjoint oriented circles
$\mu_1^\pm,\ldots , \mu_g^\pm$, a family $h=(h_1,\ldots ,h_n)$ of
chords joining the circles of $\M^+$ to those of $\M^-$, and the
signs $\varepsilon = (\varepsilon_1,\ldots ,\varepsilon_n)$ of the
chords respectively.
\end{defn}

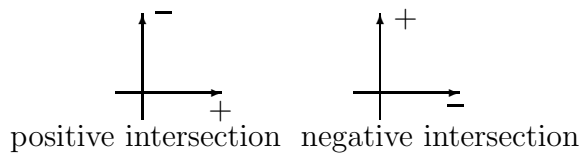
\begin{figure}[ht]
 \[\begin{picture}(200,50)
  \put(30,20){\vector(1,0){40}}
  \put(120,20){\vector(1,0){40}}
  \put(40,10){\vector(0,1){40}}
  \put(130,10){\vector(0,1){40}}
  \put(65,10){+}
  \put(45,45){{\Huge -}}
  \put(135,45){+}
  \put(155,10){{\Huge -}}
  \put(-10,0){positive intersection}
  \put(100,0){negative intersection}
 \end{picture}\]
\caption{Signs of an intersection}\label{signs}
\end{figure}

 Let $(\Sigma, M^\pm)$ be a Heegaard
diagram of a closed orientable ~3-manifold $M$. We associate to it
a Gauss diagram $G(\Sigma, M^\pm)$ in a following way: each closed
curve $m_i^\pm \in M^\pm$ can be regarded as an image of an
immersion $f^\pm\co\bigsqcup S^1\to \Sigma$. Let $\M^\pm$ be two
families of the preimages of $M^\pm$. The curves of different
families intersect, and the intersection points have two
preimages, one in each family. We join these points with a chord,
joining a circle in $\M^+$ with a circle in $\M^-$. No chords join
the circles of the same family.
 Each chord is given a sign of the
intersection, depending on whether the orientation given by the
frame (tangent to $m_i^+$ in $p$, tangent to $m_j^-$ in $p$)
coincides with the orientation of $\Sigma$, see Figure
\ref{signs}. This construction depends on the choice of the
orientation of the curves of $M^\pm$.

The choice of orientation of curves motivates the following notion
of the $\varepsilon$-equivalence:
\begin{defn}
An $\varepsilon$-move is reversing the orientation of a circle
together with reversing the signs of all chords ending on this
circle. Two Gauss diagrams $G_1$ and $G_2$ are {\em
$\varepsilon$-equivalent} if one can be obtained from the other by
a sequence of $\varepsilon$-moves.
\end{defn}
A Gauss diagram can be viewed as a picture with two rows of
circles joined with chords. For convenience, we assume the
following conventions, when drawing a Gauss diagram: the circles
of $\M^+$  are placed in the top row and are oriented
counterclockwise, while the circles of $\M^-$, placed in the
bottom row, are oriented clockwise.

\subsection{From Gauss diagrams to ribbon graphs}

Given a Gauss diagram associated to some closed orientable
~3-manifold, one might think that only the information on the
intersections is available and the information on the surface is
lost. It turns out, however, that this information is enough to
reconstruct the entire manifold. See Theorem \ref{gdthm}. This
result is based on a construction of a surface, associated to a
Gauss diagram. First we construct a ribbon graph from a Gauss
diagram, then glue discs with holes to the boundary components of
this graph.

 Let $G$ be a Gauss diagram.  Contract all chords
to points, obtaining an oriented graph $\Gamma$. Each vertex and
each edge of this graph is equipped with a sign. The signs of the
vertices are the signs of the chords they come from, where the
signs of the edges signify the family $\M^\pm$ the edge belongs
to. This graph will be the core of the ribbon graph.

\begin{figure}[ht]
\centerline{\psfig{figure=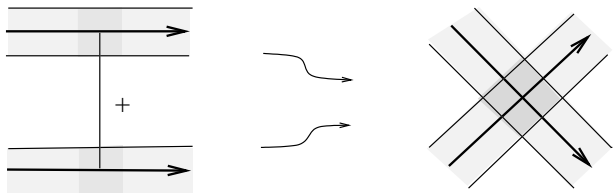,height=1in,silent=}}
\caption{Gluing the bands} \label{glue}
\end{figure}

Consider two families of $g$ oriented annuli as a product of
$\M^\pm$ with $I=[-1,1]$. Let $p^\pm$ denote the ends of a chord
$h$ in the circles $\mu^\pm$ in $\M^\pm$ respectively. Let
$U_p^\pm \subset \mu^\pm$ be small neighborhoods of $p^\pm$, their
direction induced by the orientation of $\mu^\pm$. We identify the
square $I\times I$ with $U_p^+ \times I$ by a homeomorphism
$\phi^+$ such that $\phi^+(U_p^+ \times \{0\})$ coincides with the
positive $x$ direction in $I\times I$ and the orientation is
preserved. We identify (by an orientation-preserving homeomorphism
$\phi^-$) $I\times I$ with $I \times U_p^-$ such that
$\phi^-(\{0\}\times U_p^-)$ coincides with the positive $y$
direction for a positive chord and with negative $y$ direction for
a negative chord. A homeomorphism $\phi^+\circ (\phi^-)^{-1}$
identifies the small squares in these annuli at the neighborhood
of each $p^\pm$ cross-like, as in Figure \ref{glue}.

Using such a homeomorphism for all chords, one obtains a ribbon
graph $\Gamma(G)$ with core $\Gamma$. We will call it {\em a
ribbon graph associated to the Gauss diagram}. A ribbon graph can
be viewed as a union of ribbons with bases glued to small squares,
called coupons by Reshetikhin and Turaev in \cite{rt}. We will
treat ribbon graphs in a similar way.

\begin{prop}\label{ribbon}
Let $(\Sigma, M^\pm)$ be a Heegaard diagram of some closed
orientable ~3-manifold $M$, and let $N$ be a regular neighborhood
of the union $M^+\cup M^-$ in $\Sigma$. Let $G(\Sigma, M^\pm)$ be
a Gauss diagram associated to $(\Sigma, M^\pm)$. Then
$(\Gamma(G),\Gamma)$ and $(N, M^+\cup M^-)$ are homeomorphic as
pairs.
\end{prop}
\begin{proof}We have two ribbon graphs $N$ and $\Gamma(G)$.
By construction of $\Gamma(G)$ we have the same number of coupons
and ribbons in both graphs. One can identify all the coupons
according to the order of intersections in one of the families
$M^\pm$, say $M^+$. When identifying the coupons, one has to take
care of the identification of the ribbon bases. However, the Gauss
diagram provides the exact order to do it, since the cyclic order
of the edges joining the vertex is known. Since the order of the
identification of the coupons is given by the order of edges in
$M^+$, meaning that we come to the coupon by a specific edge, the
identification of the ribbon bases is unique. Since both ribbon
graphs are oriented as surfaces and the identification goes along
annuli, this can be done in orientation-preserving way. Thus the
bases of the ribbons corresponding to the edges of $M^-$ are
identified. The cores of these ribbons are identified by the
obvious homeomorphism, hence these ribbons can be also identified
now, in the orientation-preserving way.
\end{proof}

\subsection{From ribbon graphs to splitting surfaces}
\label{cell complex}

\begin{figure}[ht]
\centerline{\psfig{figure=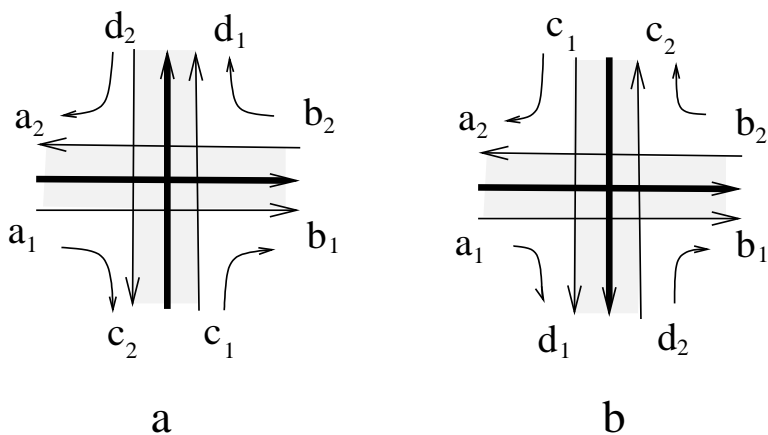,width=4in,silent=}}
\caption{Right turn rule in terms of colors} \label{colorglue}
\end{figure}

A ribbon graph $\Gamma(G)$ is an oriented surface with boundary.
As Reshetikhin and Turaev mention in \cite{rt}, we can think that
the bands of an oriented ribbon graph have two sides and we always
see only one of them. We assume that this side is such that the
boundary arc whose direction coincides with the orientation of the
edge is always seen on the right side of the edge, as in Figure
\ref{colorglue}.

If our Gauss diagram is associated to some Heegaard diagram
$(\Sigma, M^\pm)$ of a closed orientable manifold, this graph is a
ribbon graph on the surface $\Sigma$. Each family $M^\pm$
separates $\Sigma$ into a sphere with holes, hence the complement
of this graph in the surface is a collection of discs with holes.
We would like to specify the sequences of edges the boundaries of
these discs are glued along. We follow the construction of J.S.
Carter in \cite{c} for Gauss paragraphs, adapting his construction
to the case of Gauss diagrams.

\begin{figure}[ht]
\centerline{\psfig{figure=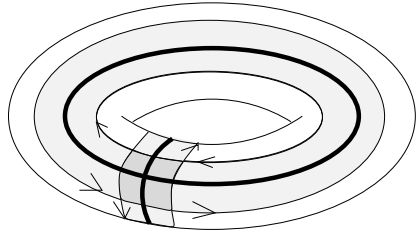,height=1.2in,silent=}}
\caption{A ribbon graph on a surface} \label{crossing}
\end{figure}

\begin{figure}[ht]
\centerline{\psfig{figure=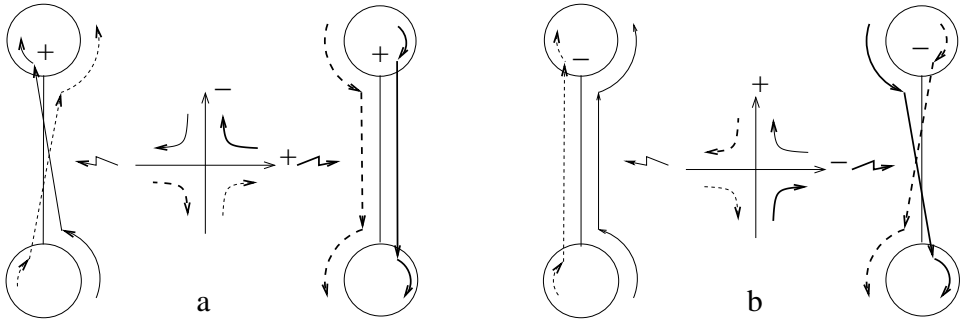,width=5in,silent=}}
\caption{Right turn rule on positive and negative intersections}
\label{negpos_fig}
\end{figure}

If one is approaching the vertex, while going along a boundary
component, he has to turn right - see Figure \ref{crossing}.  In
terms of Gauss diagrams, a turn in a vertex of $\Gamma$ means that
if one approaches a chord, going along a curve of the top family
$\M^+$, he has to descend using a chord and continue along the
curve of $\M^-$, and vice versa. This rule can be applied to any
Gauss diagram, not necessary to the one associated to some
Heegaard diagram of some manifold. The ``right turn rule" then
means that for the positive chord one has to cross it on ascending
path and to continue along it on the descending path. See Figure
\ref{negpos_fig}a. Approaching the negative chord, one ascends
along the chord and crosses it on his way downwards. See Figure
\ref{negpos_fig}b.

A sequence of edges corresponding to the boundary component is
called a {\it cycle} in a Gauss diagram. We are going to glue
discs with holes to boundary components.

Each edge of the Gauss diagram $G$ produces two arcs in the
boundary of $\Gamma(G)$. We would like to know directly from $G$
which boundary component contains each of these arcs. It proves
useful for computational purposes to include in the Gauss diagram
the information which pair of cycles the edge belongs to. An easy
way to encode it is using the elements of a finite set $C$.
Elements of $C$ are in one-to-one correspondence with boundary
components of $\Gamma(G)$.  We associate to each edge of $G$ an
ordered pair $(c_1,c_2), c_i \in C$. The first element corresponds
to the boundary component having the same direction as the edge
itself and the second one - to the opposite one.

In these terms the ``right turn rule"  will look as follows. Let
$(a_1,a_2)$, $(b_1,b_2)$ be the pairs assigned to incoming and
outgoing edges of $\mu^+$ for some chord $p$, $(c_1,c_2),
(d_1,d_2)$ - the pairs assigned to incoming and outgoing edges of
$\mu^-$ respectively. See Figure \ref{colorglue}.

 For positive intersection we have:
 $( a_1,a_2,b_1,b_2 ) =( c_2,d_2,c_1,d_1 )$, see Figure
\ref{colorglue}a,
 and for negative one we have:
 $( a_1,a_2,b_1,b_2 ) =( d_1,c_1,d_2,c_2 )$, see Figure
 \ref{colorglue}b.
From now on we incorporate this information into Gauss diagram,
without changing the notation.

We would like to construct a closed surface from $\Gamma(G)$,
hence we glue discs with holes to it. The result of this
construction is not unique. If $G = G(\Sigma, M^\pm)$, we would
like to obtain the surface homeomorphic to $\Sigma$ and to get a
Heegaard diagram producing the same manifold. To do this one has
to know which boundary components will be glued to the same disc
with holes. This leads us to the notion of a decorated Gauss
diagram.

Consider an equivalence relation on the set $C$ of all the cycles
of $G(\Sigma, M^\pm)$, where two cycles will be equivalent if they
correspond to the boundary components of the same disc with holes
in $\Sigma - N$. It is convenient to  define this equivalence
relation as a coloring of the elements of $C$ by the colors
$c_1,\ldots,c_m$ of the corresponding discs with holes.

  A {\em decorated Gauss diagram} $G^d = (\M^+,\M^-,h,\varepsilon,c)$
is a Gauss diagram with a color $c_i$  associated to each cycle,
where the colors are members of a finite set $c = (c_1,\ldots
,c_m)$. We associate to a Heegaard diagram $(\Sigma, M^\pm)$ a
decorated Gauss diagram $G^d(\Sigma, M^\pm)$. The only thing to do
is to define the coloring of the cycles. All boundary components
of the same connected component of $\Sigma - N$ have the same
color. To construct the surface associated to $G^d$, for each
color pick a sphere with a number of holes equal to the number of
cycles having this color and glue it to all these cycles. This
surface will be denoted by $S(G^d)$ and will be called {\it the
surface associated to the decorated Gauss diagram}. The manifold
associated to it can be constructed by gluing discs to $\M^\pm$
and capping the spheres with 3-balls.

\begin{prop}\label{surfdeco}
Let $(\Sigma,M^\pm)$ be a Heegaard diagram  of a closed orientable
~3-manifold $M$. The surface $S(G^d)$ associated to
$G^d(\Sigma,M^\pm)$ is homeomorphic to $\Sigma$ and the associated
manifold is homeomorphic to $M$.
\end{prop}
\begin{proof}
Since each of the families $M^\pm$ separates $\Sigma$ into discs
with holes, the surface $\Sigma -N$ consists of discs with holes,
too. By Proposition \ref{ribbon}, we have a homeomorphism
$\varphi\co\Gamma(G^d) \to N$. Consider all the cycles of
$G^d(\Sigma,M^\pm)$ colored by the same color $c_i$. All of them
correspond to the boundary components of $\Gamma(G^d)$, whose
images under $\varphi$ are glued to the same disc with holes in
$\Sigma -N$. The discs with holes glued to the boundary components
of $\Gamma(G^d)$ and the components of $\Sigma -N$ are
homeomorphic. The homeomorphism is defined by the image of
$\varphi$ on the boundaries.

 Since the Heegaard diagram was reconstructed, the manifold obtained by
 gluing discs along the core of $\Gamma(G^d)$ will be
 homeomorphic to $M$.
 \end{proof}

\begin{note}
If $\Sigma-N$ is a union of discs, then all cycles of
$G^d(\Sigma,M^\pm)$ have different colors. This means that a
manifold can be reconstructed from a non-decorated Gauss diagram,
and we can omit the decorations.
\end{note}

\subsection{From decorated Gauss diagrams back to non-decorated
ones}\label{decomit}
\begin{figure}[ht]
 \centerline{\psfig{figure=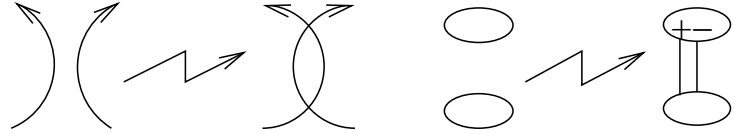,height=0.7in,silent=}}
 \caption{A pair of cancelling intersections} \label{reidmove}
 \end{figure}

 Suppose that on the surface $\Sigma$ there are two curves
 belonging to different families $M^\pm$, which can be isotoped along $\Sigma$,
 such that a pair of cancelling
 intersections appears, in a manner similar to the second
 Reidemeister move in knots. See Figure \ref{reidmove}.
 In a Gauss
 diagram this will look like adding a pair of adjacent
 chords with opposite signs.

\begin{figure}[ht]
 \centerline{\psfig{figure=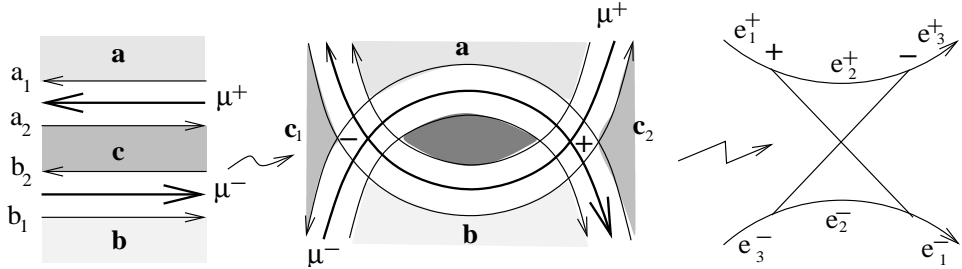,width=5in,silent=}}
 \caption{Two bands meeting with $a_2=b_2$.} \label{reidglue}
 \end{figure}

Consider a decorated Gauss diagram $G^d$. Let $e^+\in\M^+,
e^-\in\M^-$ be two edges, belonging to two cycles of the same
color. Let $(a_1,a_2),(b_1,b_2)$ be the respective colorings of
$e^+,e^-$, and let $p^+,m^+$ be points of $e^+$, $p^-,m^-$ be
points of $e^-$. The chord joining $p^+$ with $p^-$ will be
positive and the chord joining $m^+$ with $m^-$ will be negative.
Going along $e^\pm$ induces an order on the points of an edge. We
write $p<m$ if $p$ appears before $m$. The signs and the order of
the chords in $e^\pm$ depends on the colors of the edges:
\begin{itemize}
 \item if $a_1=b_1$ then $m^+<p^+$ and $p^-<m^-$
 \item if $a_1=b_2$ then $p^+<m^+$ and $p^-<m^-$
 \item if $a_2=b_1$ then $m^+<p^+$ and $m^-<p^-$
 \item if $a_2=b_2$ then $p^+<m^+$ and $m^-<p^-$.
\end{itemize}
This rule explains which sides of the bands meet when two chords
appear while two cancelling intersections are created. See Figure
\ref{reidglue} for an example of two bands meeting with $a_2=b_2$.

 We will call this change in a decorated Gauss diagram an {\em
 $R$-move}.

 During an $R$-move three new edges $e^\pm_1, e^\pm_2, e^\pm_3$ are created
 instead of each of $e^\pm$, so that $e^\pm_1$ have a common chord,
 and $e^\pm_3$ too. See Figure \ref{reidglue}. The colorings
 of these edges are given by  the following.

  If $e^\pm$ belonged to the same cycle of color $c$, then there
  is a partition of the set of all cycles colored by $c$ into two
  sets, colored into two new colors that did not exist in the
  color set of $G$. Two new cycles, containing $e_1^\pm$ and
  $e_3^\pm$ respectively, belong to different sets. A new cycle,
  consisting of two edges $e_2^+,e_2^-$ has an entirely new color
  $d$.

 If $e^\pm$ belonged to two different cycles of color $c$, then
 two new cycles appear. One  of them contains all the edges of the old cycles
 and two pairs $e^\pm_1, e^\pm_3$ and has the same color $c$. The other one
 consists only of $e_2^+,e_2^-$ and has a new color $d$.

 The edges $e_1^\pm$ and $e_3^\pm$ have the same coloring pairs as
 $e^\pm$, save that $c$ can be replaced by a new color, as
 mentioned above. The edges $e_2^\pm$ have a new coloring,
 obtained from the coloring of $e^\pm$ by the following changes:
 $c$ is replaced by $d$, the second color in the pairs is switched
 between $e^+$ and $e^-$, and the colors inside each pair are also
 switched. E.g., if $e^+$ is colored by $(a,c)$ and $e^-$ - by
 $(b,c)$, as in Figure \ref{reidglue}, the edge  $e_2^+$ is
 colored by $(d,b)$ and  $e_2^-$ - by $(d,a)$.

 The opposite move to an $R$-move is called {\em $R^{-1}$-move}.

 \begin{defn}
 Two decorated Gauss diagrams are called
 {\em $R$-equivalent} if one can be obtained from the other by a
 sequence of $R$ and $R^{-1}$-moves.
 \end{defn}

 \begin{note}
 $R$-equivalent diagrams have the same Heegaard splitting.
 \end{note}
 \begin{note}
 An $R$-move reminds us the second Reidemeister move in knots. There is
 no equivalent to the third Reidemeister move, since three curves
 intersecting would produce an intersection between the curves of
 the same family, which is forbidden.\pole
 \end{note}

 A ribbon graph $\Gamma(G^d)$ may be connected,
 and may be not. This leads us to the notion of $R$-connected
 diagram. Let $G^d$ be an abstract decorated Gauss diagram.
 Consider a graph, whose vertices are connected components of
 $G^d$. Two vertices are connected with an edge if and only if the
 sets of colors appearing on these connected components have a
 color in common. We call the diagram $G^d$ {\em $R$-connected} if and
 only if this graph is connected.

 \begin{lem}
 The surface $S(G^d)$ is connected if and only if $G^d$ is
 $R$-connected.
 \end{lem}
 \begin{proof}
 The connected components of $\Gamma(G^d)$ are exactly the
 connected components of $G^d$ itself. Two components will have
 the same disc with holes glued to both of them if and only if
 they have a color in common. Hence the connected components of
 the graph defined above correspond to connected components of $S(G^d)$.
 \end{proof}

\begin{figure}[ht]
\centerline{\psfig{figure=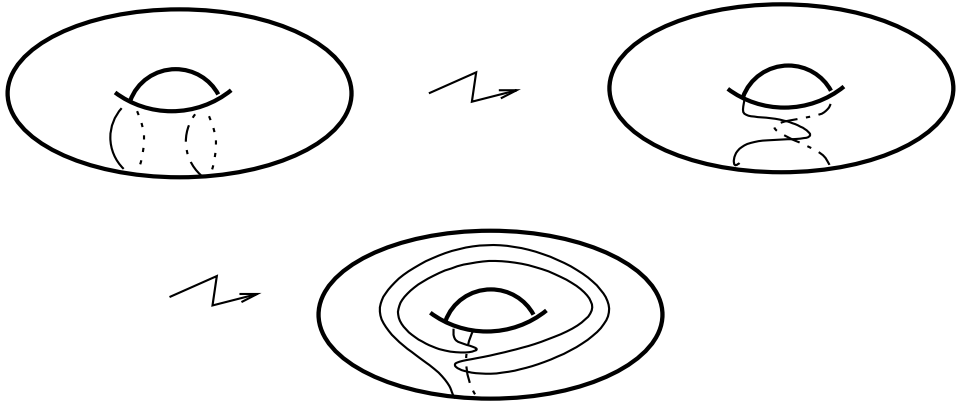,width=5in,silent=}}
\caption{Adding pairs of cancelling intersections changing annuli
into discs} \label{reidemeister}
\end{figure}

 \begin{prop}\label{requiv}
  Any decorated Gauss diagram $G^d$ is $R$-equivalent to a
 Gauss diagram with the number of colors used on it equal to the number of
 cycles. The resulting diagram is connected if and only if $G^d$
 is $R$-connected.
 \end{prop}
 \begin{proof}
 Consider a decorated Gauss diagram $G^d$.
 Let $\Sigma_i$ be a non-simply connected component of $S(G^d)-\Gamma(G^d)$.
 It is a disc with holes and each of its boundary components contains arcs of both
 families.
 Perform an $R$-move between two different boundary components
 of $\Sigma_i$. Such a
 move produces two components, such that one of them is a disc
 with a new color. The other one has one hole less than the
 initial, thus lowering the number of cycles colored by
 the same color. See Figure \ref{reidemeister}. An $R$-move does not change
 the number of boundary components of any
 other component of $S(G^d)-\Gamma(G^d)$.
  The resulting diagram will be
 $R$-equivalent to the previous one. One can continue this process
 till $\Sigma_i$ turns into a disc.

 Remove all non-simply-connected components of $S(G^d)-\Gamma(G^d)$
 by the $R$-moves. This will give a decorated Gauss
 diagram with all cycles colored differently.

 If the Gauss diagram $G^d(\Sigma,M^\pm)$ is not connected, we
assume for convenience that it has two connected components. If
$G^d$ is $R$-connected, there are two cycles colored by the same
color in different components of $G^d(\Sigma,M^\pm)$. An $R$-move
between these cycles will change the diagram into a connected one.
 \end{proof}

\begin{cor}\label{connexist}
For any Heegaard diagram $(\Sigma,M^\pm)$ of a closed connected
orientable ~3-manifold there is a connected Gauss diagram $G$ with
all cycles colored in different colors $R$-equivalent to
$G^d(\Sigma,M^\pm)$. Moreover, there is a Heegaard diagram
isotopic to $(\Sigma,M^\pm)$, such that $G$ is associated to it.
 \end{cor}

\subsection{Gauss diagrams and their associated surfaces}\label{surf_bdry}
To compute the genus of the surface associated to an abstract
Gauss diagram, we will use the notion of the {\em excess of a
coloring}, which is the difference $\Delta_c$ between the number
of cycles and number of colors. Obviously, if all cycles have
different colors, $\Delta_c$ vanishes.
\begin{prop}\label{genus}
The genus of a surface $S(G^d)$ associated to an abstract
decorated Gauss diagram $G^d$ is $g_{S(G^d)} = 1+ \Delta_c +
\frac{1}{2}(|h|-|c|)$, where $|c|$ is the number of cycles and
$|h|$ is the number of chords in $G^d$.
\end{prop}
\begin{proof}
 This can be done by computing the Euler characteristic of $S(G^d)$. The
 surface $S(G^d)$
 is constructed by gluing discs with holes to the boundary components of a ribbon graph
 $\Gamma(G^d)$. The number of vertices of $\Gamma(G^d)$
 is the number of chords $|h|$. The number of edges of $\Gamma(G^d)$ equals twice
 the number of chords. We do not count circles that have no chord attached to
 them, since they produce annuli in $\Gamma(G^d)$, and their Euler
 characteristic is $0$. Thus $\chi(\Gamma(G^d))= |h|-2|h|=-|h|$.
 Glue discs with holes to $\Gamma(G^d)$. For each of them, the
 number of the boundary components equals the number $|c_i|$ of
 cycles colored by $c_i$. The total number of these discs with
 holes is the number of colors $|c|-\Delta_c$. The Euler
 characteristic of each such disc with holes is $\chi_i =
 2-|c_i|$. Summarizing, we have:
 \[ \chi(S(G^d)) = -|h|+ \sum^{|c|-\Delta_c}_{i=1} \chi_i =
 -|h|+2|c|-2\Delta_c  + \sum^{|c|-\Delta_c}_{i=1} |c_i|.\]
 The sum over $|c_i|$ equals the number of cycles $|c|$, giving
 $\chi(S(G^d))= -|h|+|c|-2\Delta_c$.
 The genus is given by $\chi(S(G^d)) = 2-2g_{S(G^d)}$.
 \end{proof}
 A surface associated to a
non-decorated Gauss diagram $G$ is a surface obtained by gluing
discs to all boundary components of $\Gamma(G)$. In other words,
this is the surface associated to a decorated Gauss diagram with
all cycles colored in different colors.
\begin{cor}\label{gnondeco}
The genus of the surface associated to the connected non-decorated
Gauss diagram is  $g_{S(G)} = 1+ \frac{1}{2}(|h|-|c|)$.
\end{cor}
\begin{prop}
 Let $(\Sigma,M^\pm)$ be a Heegaard diagram  of some
 closed orientable ~3-manifold $M$. Assume that $G(\Sigma,M^\pm)$ is
 connected. Let $g$ be the number of curves in each family
 $M^\pm$.
 If $g_{S(G)}= g$, then $\Sigma-(M^+\cup M^-)$ is a collection of
 discs, otherwise
 $(\Sigma,M^\pm)$ is reducible.
\end{prop}
\begin{proof}
 If $g_{S(G)}= g$, the surface associated to a non-decorated Gauss
 diagram is homeomorphic to $\Sigma$. The boundaries of the meridional
 disc system are
 also reconstructed, hence $\Sigma$ is separated into discs by the curves of
 $(M^+\cup M^-)$.

 If $g_{S(G)}\neq g$, then cutting $S(G)$ along the curves of the diagram
 produced surfaces
 other than discs. Thus any essential closed curve in such a non-simply
 connected component bounds a disc
 in both handlebodies of the Heegaard splitting. This means that the
 Heegaard splitting is reducible.
\end{proof}

\subsection{Equivalent diagrams}

We introduce now two types of moves of Gauss diagram motivated by
some equivalences of Heegaard splittings.

\begin{figure}[ht]
\centerline{\psfig{figure=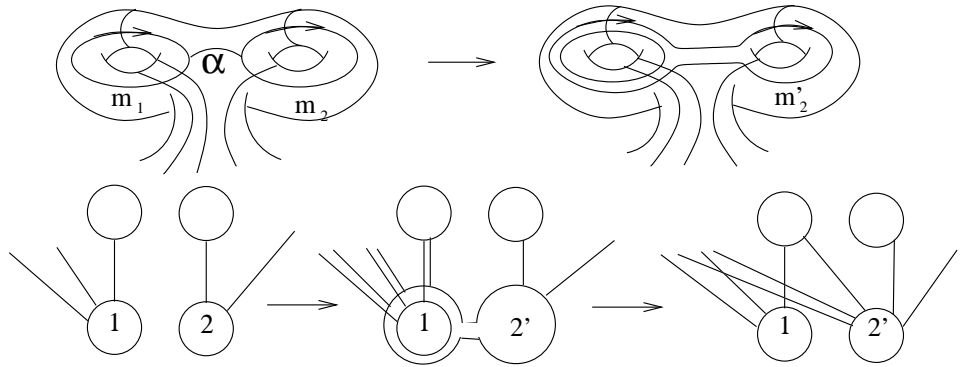,width=5in,silent=}}
\caption{Handle slide} \label{handlesliding}
\end{figure}

Suppose that two curves $m_1,m_2$ on the Heegaard diagram
$(\Sigma,M^\pm)$, belonging to the same family $M^\pm$, say $M^+$,
can be joined together with a simple curve $\alpha$ disjoint from
all other curves of $M^\pm$. Let $\alpha',\alpha''$ be the
boundary components of a thin collar $N(\alpha)=\alpha \times
[-1,1]$ of $\alpha$ in $\Sigma-(M^+\cup M^-)$. Set
\[m_2' =
(m_2 - N(\alpha))\cup (m_1 - N(\alpha))\cup \alpha'\cup \alpha''.
\]
 By a slight push make this curve
disjoint from $m_1$. This does not change the handlebodies of the
splittings, only the meridional discs defining them. The Heegaard
diagram $(\Sigma,(M^+-m_2)\cup m_2', M^-)$ is obtained from the
diagram $(\Sigma,M^\pm)$ by a handle slide of $m_2$ along $m_1$.
See Figure \ref{handlesliding}.

This process changes the Gauss diagram of $M$, with a small
section of $\mu_2$ being substituted by a curve equal to $\mu_1$
without a small segment. All the curves that intersected $m_1$ now
intersect $m'_2$ also, which means that we have the copies of all
chords that joined $\mu_1$ now joining $\mu'_2$. The signs of the
intersections (and the chords) are the same, if the direction of
the segment removed coincides with the direction of the added
segment. If not, one has to change the direction of the added
segment, and reverse the signs for all chords joining it, as in
the case of $\varepsilon$-equivalence. Clearly, the move is
allowed only if the segments removed from both cycles are of the
same color.  We say that the resulting diagram is obtained from
the initial one by an {\em $H$-move}. The opposite of this move
will be called {\em $H^{-1}$-move}.

\begin{figure}[ht]
\centerline{\psfig{figure=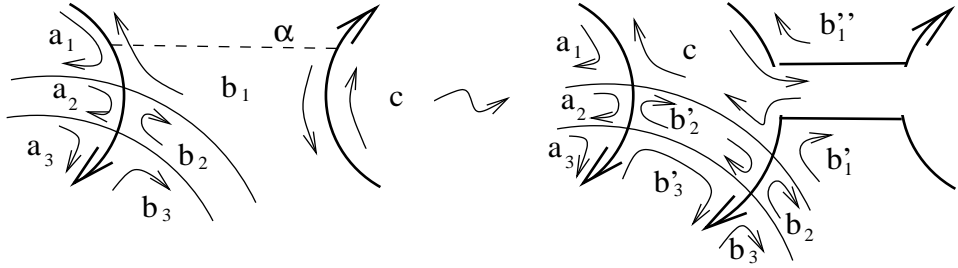,width=5in,silent=}}
\caption{Color change during the handle slide} \label{colorslide}
\end{figure}

The change of the colors in the case of handle slide is
illustrated in Figure \ref{colorslide}. The color of the region is
the color of the cycle bounding it. Let $(a_i,b_i)$ be the
colorings of the edges $e_i$ of $\mu_1$, and let $\alpha$ join
$e_1$ and $e\in \mu_2$. Let $(b_1,c)$ be the coloring of $e$. The
edges $e_i$, save $e_1$, will now be parts of new small cycles,
consisting only of $4$ edges each and having new colors $b'_i$.
The edges $e_i, i>1$ will have colors $(a_i,b'_i)$. The edge $e_1$
will have the coloring $(a_1,c)$, if the directions of the added
and removed segments coincide. In this case also the edges $e'_i$
added to $\mu_2$ have colors $(b'_i,b_i)$ and belong to the same
new cycles. The edges $e'_i$ have colors $(b_i,d_i)$. The
remaining two edges of these new cycles have the colors
$(b'_i,b'_{i+1})$ or $(b'_{i+1},b'_i)$, depending on the signs of
their ends. We described already the colors of all $e'_i$ save the
first and the last ones. Their coloring equals that of $e$ if
$\alpha$ joined two edges of different cycles. If not, the color
$b_1$ splits and the coloring is as in Figure \ref{colorslide}. If
the direction of the added segment is opposite to that of the
removed one, the colorings of the edges of the added segment
switch.
\begin{lem}\label{hgs}
$H$-move does not change the genus of $S(G^d)$.
\end{lem}
\begin{proof}
Consider the slide of $\mu_i$ along $\mu_j$. Let $|h_j|$ be the
number of chords joining $\mu_j$. Then the total number of chords
becomes $|h|+|h_j|$. The number of cycles and colors depend on
whether the edges joined by $\alpha$ belonged to the same cycle or
not. If $\alpha$ connected two edges of the same cycle, two cycles
are created instead of it. Their colors are different, and the set
of the cycles colored by the initial color is partitioned into two
sets, each with a new color. Additional $(|h_j|-1)$ new $4$-edged
cycles are also created, each with a separate color. Hence both
the number of colors and the number of cycles increases by
$|h_j|$, and $\Delta_c$ does not change. The formula of
Proposition \ref{genus} for the changed diagram ${G'}^d$ gives:
$g_{S({G'}^d)}=
1+\Delta_c+\frac{1}{2}(|h|+|h_j|-|c|-|h_j|)=g_{S(G^d)}$. See also
Figure \ref{colorslide}.

If $\alpha$ connected two edges of two different cycles of the
same color, these cycles are now replaced by a single one. Hence
the number of colors increases only by $(|h_j|-1)$, while the
number of cycles increases by $(|h_j|-2)$. The formula of the
Proposition \ref{genus} then gives us:  $g_{S({G'}^d)}=
1+\Delta_c-1+\frac{1}{2}(|h|+|h_j|-|c|-|h_j|+2)=g_{S(G^d)}$.
\end{proof}

\begin{figure}[ht]
\centerline{\psfig{figure=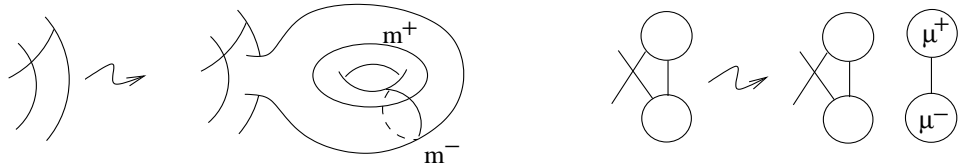,width=5in,silent=}}
\caption{The stabilization on the Heegaard and Gauss diagrams}
\label{stabilization}
\end{figure}

Another basic operation with Heegaard splittings is stabilization.
This operation consists of adding a pair of cancelling handles.
See Figure \ref{stabilization}. One can think of it as a connected
sum of $M$ with $S^3$, where the ball of connection is a small
neighborhood of a point on the splitting surface $\Sigma$ disjoint
from the curves of $M^\pm$, and $S^3$ is thought as having a
standard genus 1 splitting $T^\pm$. The boundaries of the
meridional discs of the solid tori $T^\pm$ intersect in one point
and are disjoint from $M^\pm$. This means that in terms of Gauss
diagrams the stabilization can be viewed as adding a circle to
each family $\M^\pm$. These two circles are joined with one chord
(the sign depends on the orientation of the boundaries of the
meridional discs of $T^\pm$ and can be chosen arbitrarily), and no
chord joins them to the other circles. They produce a single cycle
whose color is one of the existing colors. The color represents
the place where the stabilization was made. We call this operation
an {\em $S$-move}. There is, certainly, the inverse operation of
removing two cancelling handles, disjoint from all others. This
looks like removing a pair of circles joined by a single chord and
not connected to any other circle of a Gauss diagram. We call it
{\em $S^{-1}$-move}.

\begin{lem}\label{sgs}
 $S$-move increases the genus of $S(G^d)$ by $1$.
\end{lem}
\begin{proof}
Let ${G'}^d$ be a stabilization of $G^d$.  A single cycle was
added to $G^d$, and no new colors. Hence $\Delta_c$ increased by
$1$. We also have one chord more. Substituting to the formula of
the Proposition \ref{genus} for $g_{S({G'}^d)}$ we have:
$g_{S({G'}^d)}=1+\Delta_c+1-\frac{1}{2}(|h|+1-|c|-1)=g_{S(G^d)}+1$.
\end{proof}
\begin{defn}
Two Gauss diagrams are called {\em $H\hspace{-2pt}S$-equivalent}
if one can be obtained from the other by a sequence of $S$-moves,
$H$-moves and their opposites.
\end{defn}

\subsection{From splitting surfaces to ~3-manifolds}

One of the problems of Heegaard splittings is that a ~3-manifold
can be presented by different Heegaard splittings (possibly of a
different genus). Moreover, the same Heegaard splitting can be
presented by different Heegaard diagrams. The
$H\hspace{-2pt}S$-equivalence reflects this phenomenon on the
level of Gauss diagrams.

\begin{defn}
Two Gauss diagrams are {\em equivalent} if one can be obtained
from the other by a sequence of $\varepsilon,R,H,S$-moves or their
opposites.
\end{defn}

\begin{them}\label{gdthm}
A closed connected orientable ~3-manifold is in one-to-one
correspondence with the equivalence class of its Gauss diagram.
\end{them}
\begin{proof}
 A well known theorem by Singer says that one can pass from one
 Heegaard splitting to any other Heegaard splitting of the same
 closed manifold by a consequent use of handle slides and stabilizations.
 See \cite{s}. On the level of Gauss diagrams it means that the
 Gauss diagrams associated to these Heegaard splittings are
 $H\hspace{-2pt}S$-equivalent.

 Now let $(\Sigma,M^\pm)$ be a Heegaard diagram of some closed orientable
 ~3-manifold M.
 By Corollary \ref{connexist}, there is a connected Gauss diagram
 $G$ with all cycled colored differently, $R$-equivalent to
 $G(\Sigma,M^\pm)$, and the curves of the Heegaard diagram can be isotoped
 so that $G$ is the diagram associated to $(\Sigma,M^\pm)$ after
 the changes. We will not change the notation for the isotoped
 $M^\pm$.
 Construct the surface associated to $G$. By Proposition
 \ref{ribbon}, a ribbon graph associated to $G(\Sigma,M^\pm)$
 is homeomorphic to a regular neighborhood of $M^\pm$. Consider
 the boundary components of this neighborhood. Each of them is a simple closed curve
 disjoint from the curves of $M^\pm$, homotopic to the image of a cycle
 in $G(\Sigma,M^\pm)$ under this homeomorphism, and all cycles of
 $G(\Sigma,M^\pm)$ are presented by these boundary components. If this curve is
 essential in $\Sigma$, compress along it. Continue for all cycles. Take
 only those components of the resulting surface that contain
 curves of $M^\pm$. This surface is homeomorphic to the
 surface $S(G)$.

 The complement of the ribbon graph in $S(G)$ is a
 collection of discs. In the surgered surface containing
 $M^\pm$, the complement of $M^\pm$ is also a collection of
 discs, the boundary of each is a cycle in the Gauss diagram. A
 homeomorphism is constructed by taking the corresponding discs to
 each other and the correspondence is determined by the boundary.
 Thus any two Heegaard diagrams with equivalent Gauss diagrams can
 be obtained from a surface associated to their common Gauss
 diagram by a connected sum with unknotted tori and handle slides.
\end{proof}

\subsection{Visualizing Gauss diagrams}\label{GDvisualizing}
The algorithm of Section \ref{cell complex} gives us the
possibility to reconstruct the surface from the Gauss diagram, but
it is hard to visualize it. The following algorithm allows us to
reconstruct the picture of the Heegaard diagram in an easy and
quite visual way, provided that a Gauss diagram comes from some
Heegaard diagram of a closed orientable ~3-manifold M and is
connected. We also need to assume that $g=g_{S(G)}$, where
$g_{S(G)}$ is given by the formula of Corollary \ref{gnondeco}.

Consider such a Gauss diagram $G$. Since one of the families
$M^\pm$ of curves on the surface can be taken to be standard, we
start with an oriented disc with $2g-1$ holes. Let the boundary
have the orientation induced by the orientation of the disc. Split
the set of the boundary components in pairs. Choose an
orientation-reversing homeomorphism between the curves and
identify them in pairs.
\begin{figure}[ht]
\centerline{\psfig{figure=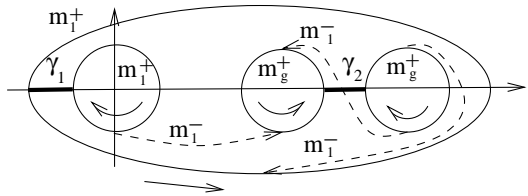,height=1in,silent=}}
\caption{Surface $S$ cut along the curves of $M^+$} \label{gdv}
\end{figure}

This gives us a surface $S$ of genus $g$. See Figure \ref{gdv},
where the boundary components already have reversed orientations,
where necessary. After the identification we refer to the
identified curves as $m^+_i, \ i = 1,\ldots,g$.

 The curves $m^+_i$ can be considered as a map from a Gauss diagram to a
 Heegaard diagram, defining the closed curves belonging to one of
 the families. To get a map from $\M^-$ to $S$ for all segments of the
 circles of  $\M^-$ between the endpoints of the chords draw
 disjoint curves in $S$ between the points corresponding to the
 upper ends of the chords intersecting this segment. When choosing the copy of
 $m^+_i$ that should contain the final point of the edge of $m^-_j$, one has
 to take the chord sign into account. For instance, taking the orientations of
 $m_i$'s
 as in Figure \ref{gdv}, $m^-_1$ has to join the counterclockwise copy of
 $m^+_2$ for a negative intersection and the clockwise copy of
 $m^+_2$ for a positive one. The
 concatenation of these curves after identifying the $m_i$'s
 will be a closed curve in $S$ corresponding to a circle of $\M^-$.

The Gauss diagram associated to the Heegaard diagram obtained
above is the same as $G$. By Theorem \ref{gdthm}, the manifold
given by this Heegaard diagram is homeomorphic to $M$.

\subsection{First simple computations using Gauss diagrams}\label{comput1}
The fundamental group of a closed orientable ~3-manifold $M$ and
its first homology group can be computed easily from a Gauss
diagram $G$, associated to some Heegaard splitting of $M$. For
each circle $\mu^-_j$ of $\M^-$ we write a word $r_j$ in letters
$g_i, i = 1,\ldots,g$ in the following way: let
$h_1^{\varepsilon_1}, \ldots, h_n^{\varepsilon_n} $ be the
sequence of chords with signs $\varepsilon_i \in \{ \pm 1 \}$
joining $\mu^-_j$ to the circles of $\M^+$. Suppose a chord
$h_i^{\varepsilon_i}$ joins $\mu^-_j\in \M^-$ to $\mu^+_k \in
\M^+$. For such a chord we write $g_k^{\varepsilon_i}$ in $r_j$.
Thus we get a cyclic word, which is a relation in the presentation
of $\pi_1(M)$.

\begin{prop}\label{pi1}
 Let $(\Sigma,M^\pm)$ be a Heegaard diagram of some closed orientable
 ~3-manifold $M$. Assume that $G= G(\Sigma,M^\pm)$ is connected and
 $g=g_{S(G)}$. Then $ \pi_1(M)$ has a presentation $\langle
g_1,\ldots,g_g|r_1,\ldots,r_g \rangle$,
 where $g_i$ correspond to the
circles of $\M^+$ and $r_j$ correspond to the circles of $\M^-$
and are obtained as above.
\end{prop}
\begin{proof}
 We use the algorithm and the notation of Section \ref{GDvisualizing}.
Each $g_i$ corresponds to a curve $\gamma_i$, see Figure
\ref{gdv}. The fundamental group $\pi_1(M)$ can be obtained from
the fundamental group of $S$ by adding the relations given by the
discs glued to the curves of the Heegaard diagram. The curves
$m^+_i,\gamma_i,\ i=1,\ldots,g$ form the generator set for
$\pi_1(S)$. When we glue discs to the curves corresponding to
$\M^+$, we remain with $\gamma_i$'s only. They will be the
generators of $\pi_1(M)$. When one wants to count how many times
does the curve wind along the handle defined by $\gamma_i$, he
simply counts how many times it crosses some fixed curve dual to
$\gamma_i$, which is $m^+_i$. This explains the form of the
relations.
\end{proof}
Let $a_{ij}$ denote the algebraic number of chords joining the
circle $\mu^-_j$ in $\M^-$ to a circle $\mu^+_i$ in $\M^+$.
\begin{cor}\label{H1(M)}
The first homology group of $M$ is the abelian group generated by
$g_1,\ldots,g_g$, where $g_i$ correspond to the circles of $\M^+$.
The relations  are $r_1,\ldots,r_g$, where $r_j =\sum_i a_{ij}g_i
= 0$.
\end{cor}
The algebraic numbers of chords joining the circles can also be
written in the intersection matrix $A_{ij}$. This matrix allows us
to determine whether a manifold is a homology sphere.

\subsection{Examples}\label{ex-clmfld}
\begin{description}
\item[The 3-sphere] $S^3$ is the simplest example. We represent a Gauss
 diagram of the standard genus $1$ splitting of it. See Figure
 \ref{sphere}.

 \begin{figure}[ht]
 \centerline{\psfig{figure=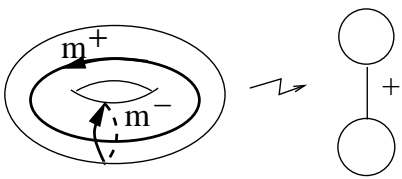,height=0.9in,silent=}}
 \caption{Heegaard and Gauss diagrams of $S^3$} \label{sphere}
 \end{figure}
\item[Lens spaces] They are a bit more complicated. In Figure \ref{lens}
 we present the Gauss diagrams of $L(5,1)$ and $L(5,2)$. The signs
 of all chords are $+$ on both diagrams. The
 fundamental group is the same, as one should expect:
 \[\pi_1(L(5,q)) = \langle g|g^5\rangle = Z_5, \ q=1,2.  \]
 \begin{figure}[ht]
 \centerline{\psfig{figure=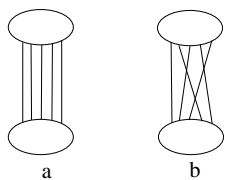,height=0.9in,silent=}}
 \caption{Gauss diagrams of $L(5,1)$(a) and $L(5,2)$(b)} \label{lens}
 \end{figure}
\item[Poincare homology sphere] Figure \ref{poincare} represents a
Heegaard diagram and the associated Gauss diagram of the Poincare
manifold $P^3$, as it appears in \cite[p.245]{r}. The picture
presented by Rolfsen looks exactly like the one obtained by the
algorithm of Section \ref{GDvisualizing} from a Gauss diagram. The
fundamental group, according to the Proposition \ref{pi1}, is
\[\langle g_1,g_2|
    g_1^{-4} g_2 g_1 g_2, g_1g_2^{-2}g_1g_2\rangle = \langle g_1,g_2|(g_1g_2)^2
    = g_1^5 = g_2^3\rangle,
\]
which is the binary icosahedral group, as expected. Abelinizing,
we get $H_1(P^3) = \langle e\rangle$, hence $P^3$ is indeed
homology sphere.

\begin{figure}[ht]
\centerline{\psfig{figure=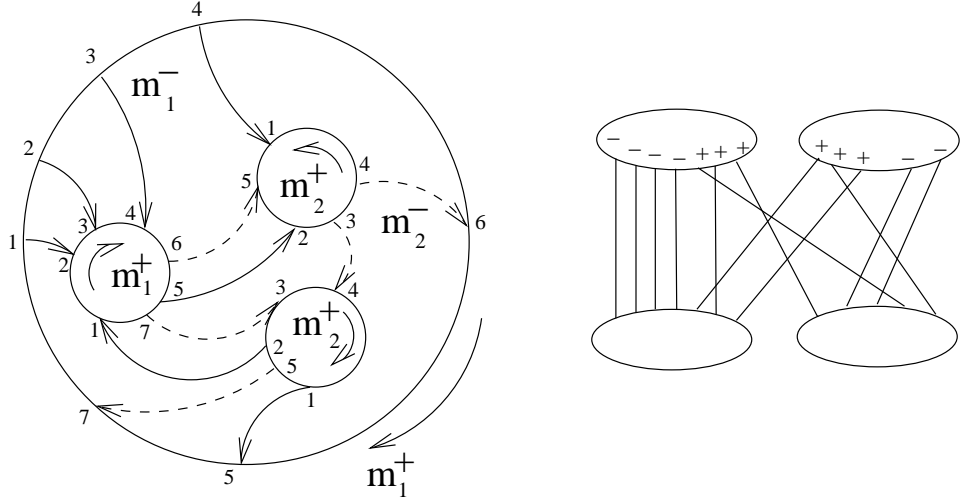,width=5in,silent=}}
\caption{Heegaard and Gauss diagrams of $P^3$} \label{poincare}
\end{figure}

\item[Another homology sphere]We consider the example of
 a homology sphere given by Hempel in \cite[p.19]{he}. Its
 Heegaard diagram is shown in Figure
 \ref{homsph}, together with its Gauss diagram. One can easily
 calculate the fundamental group:
 \[ \pi_1(M) = \langle g_1,g_2|
    g_1g_2^{-1}g_1^{-1}g_2^2g_1^{-1}g_2^{-1},
    g_1g_2g_1g_2^{-1}g_1^{-1}g_2^{-1}\rangle ,  \]
 which gives trivial first homology group. The intersection matrix
 is in this case
 \[ \left(
 \begin{array}{cccc}
 0&0&-1&1  \\
 0&0&0&-1  \\
 -1&0&0&0  \\
 1&-1&0&0  \\
 \end{array}
 \right), \]
and its determinant is $1$, in agreement with the expectations.
\begin{figure}[ht]
 \centerline{\psfig{figure=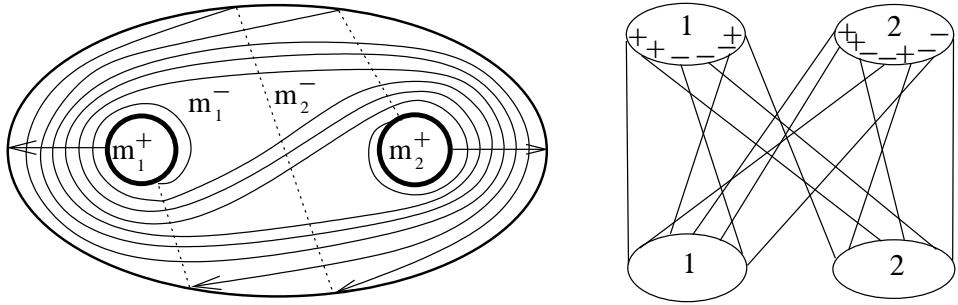,width=5in,silent=}}
 \caption{Heegaard and Gauss diagrams of a homology sphere from \cite{he}}
  \label{homsph}
 \end{figure}

\end{description}

\section{Abstract Gauss diagrams}\label{abstgd}
A natural question that arises when dealing with Gauss diagrams of
~3-manifolds is: given an abstract Gauss diagram, can one find a
manifold that this Gauss diagram is associated to it?

Gauss diagrams for ~3-manifolds closely resemble Gauss diagrams
for knots and links. The equivalence classes of Gauss diagrams for
links that do not necessarily correspond to a link are called
virtual links. We can also call the equivalence classes of Gauss
diagrams of ~3-manifolds virtual manifolds.


Below we show that a ~3-manifold, in general with a boundary, can
be associated to any abstract Gauss diagram. In Section \ref{bdry}
we give the conditions for an abstract Gauss diagram to represent
a closed ~3-manifold or a component of a knot in some closed
orientable ~3-manifold.

\subsection{Gauss diagrams for ~3-manifolds with boundary}

In order to generalize the definition of Heegaard splitting to
~3-manifolds with boundary, one replaces the handlebodies with the
compression bodies. See, e.g. \cite{sc}. A {\em compression body}
$H$ is a ~3-manifold with boundary obtained in a following way:
let $\Sigma$ be a closed surface. Take a cylinder $\Sigma\times
[-1,1]$, attach to $\Sigma\times \{-1\}$ a collection of
2-handles. If the resulting manifold has any spherical components
of the boundary, fill them with 3-balls. The surface $\Sigma\times
\{1\}$ is called the outer boundary $\partial^+ H$ of the
compression body, and the surface $\partial^- H = \partial H -
\partial^+ H$ is called the inner boundary of $H$. If a boundary
has one component only, namely its outer boundary $\partial^+ H$ ,
we get a handlebody. The boundary of a compression body may
consist of more than two components. The genus of the outer
boundary is the genus of a compression body. Unlike the case of a
handlebody, the genus does not determine the manifold completely.
Two compression bodies of the same genus have homeomorphic outer
boundaries, but there are no restrictions on the inner boundary,
and the manifolds can be completely different. A Heegaard
splitting of a manifold with boundary is then a union of two
compression bodies $H^\pm$ glued using a homeomorphism of the
outer boundaries $\partial^+ H^\pm$ of the compression bodies. The
surface $\partial^+ H^+ = \partial^+ H^-$ is then called a
splitting surface of the Heegaard splitting. We also allow now the
families $\M^\pm$ to have different number of circles. The number
of circles in each family $\M^\pm$ will be denoted by $g^\pm$.

Assume that $G^d$ is an abstract decorated Gauss diagram with $g$
circles in each family. Then one can look at the Gauss diagram as
representing curves on $S(G^d)$ to which the discs are glued.
Being more precise, one takes a collar $S(G^d)\times [-1,1]$ and
attaches 2-handles to it along the curves in $\M^+\times\{1\}$ and
$\M^-\times \{-1\}$. If a 3-sphere appears in the boundary, it is
filled with a 3-ball. Then $S(G^d)\times \{ 0\}$ is a Heegaard
splitting surface for a manifold, possibly with boundary.

If one adds to one of the families of a Heegaard diagram a curve
contractible on the splitting surface, gluing disc along it
produces a spheric boundary component. This sphere is filled with
a ball and the Heegaard splitting does not change. This move is
called a bubble-move or a {\em $B$-move}. In terms of Gauss
diagrams this looks like the addition to one of the families
$\M^\pm$ a circle (called a {\em bubble}) without chords. This
circle produces two cycles, one colored in a new color and the
other in one of the existing colors. The choice of the existing
color is the choice of the component of $S(G^d)-\M^\pm$ to which
the bubble is added. The opposite move to this one, i.e. a
removing of such a circle, is called $B^{-1}$-move.

\begin{defn}
Two abstract Gauss diagrams are called {\em $B$-equivalent} if one
can be obtained from the other by a sequence of $B$ and
$B^{-1}$-moves.
\end{defn}

The change on the Gauss diagram produced by a $B$-move is less
significant than that made by an $\varepsilon$-move. Moreover, a
$B$-move allows us to get rid of $\varepsilon$-move.

\begin{figure}[ht]
 \centerline{\psfig{figure=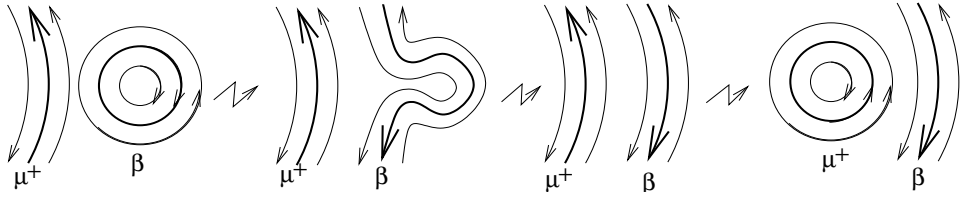,width=5in,silent=}}
 \caption{$\varepsilon$-move as bubble moves and handle slides}
  \label{epsilon}
 \end{figure}

\begin{prop} The $\varepsilon$-equivalence relation is generated
by the $H\hspace{-2pt}S$ and $B$ equivalence relations.
\end{prop}
\begin{proof}
Let $G^d_1,G^d_2$ be two abstract decorated Gauss diagrams. Assume
that $G^d_2$ can be obtained from $G^d_1$ by an
$\varepsilon$-move. We would like to show that $G^d_2$ can be
obtained from $G^d_1$ by a sequence of $B,B^{-1}$ and $H,
H^{-1}$-moves. The diagram $G^d_1$ differs from $G^d_2$ by the
signs of the chords joining one circle, say $\mu^+$ in $\M_1^+$.
Choose an edge in $\mu^+$. Let $(c,d)$ be its coloring. Perform a
$B$-move adding a circle $\beta$ to $\M_1^+$, colored by $(e,c)$,
where $e$ is a new color. Perform an $H$-move, sliding $\beta$
along $\mu^+$. The new circle will have the direction opposite to
that of $\mu^+$, and the chords will have opposite signs. Perform
an $H^{-1}$-move on $\mu^+$ changing it into bubble. Remove this
bubble with $B^{-1}$-move. The resulting diagram is $G^d_2$. See
Figure \ref{epsilon}.
\end{proof}

\begin{lem}\label{bgs}
$B$-moves do not change the genus of $S(G^d)$.
\end{lem}
\begin{proof}
Let ${G'}^d$ denote the diagram resulting from a Gauss diagram
$G^d$ after the $B$ move. A $B$-move adds a circle to one of the
families, without any chords joining it. Two new cycles and one
new color are added, increasing the excess of coloring $\Delta_c$
by $1$ and $|c|$ by $2$. Hence
$g_{S({G'}^d)}=1+\Delta_c+1+\frac{1}{2}(|h|-|c|-2)=g_{S(G^d)}$.
\end{proof}

\subsection{The boundary and its genus}\label{bdry}
  We would like to understand which conditions a Gauss
diagram should satisfy so that a manifold associated to it is
closed. If the number of circles in each of $\M^\pm$ equals
$g_{S(G^d)}$, we have an honest Heegaard splitting of a closed
orientable ~3-manifold, provided that neither $\M^+$ nor $\M^-$
separate the surface. See, e.g., \cite{fm}.

For a decorated Gauss diagram $G^d$ consider a pair of graphs
$C^\pm=C^\pm(G^d)$. The vertices of both graphs are the colors
used in $G^d$. Two vertices $p,q$ of, say, $C^+$ have a common
edge if there is an edge or a circle without chords in $\M^-$ with
a coloring $(p,q)$ or $(q,p)$. The following Lemma is obvious:
\begin{lem}\label{separation}
The connected components of $C^\pm$ correspond to the connected
components of $S(G^d)-\M^\pm$ respectively.
\end{lem}

Given a Gauss diagram $G^d$, we will denote by $\pg^\pm$ the sum
of the genera of the components of $\partial^-H^\pm$ respectively.
\begin{prop}\label{pg}
Let $G^d$ be a Gauss diagram, and let $k^\pm$ be the number of
components in $C^\pm(G^d)$. Then $\pg^\pm = k^\pm -g^\pm
+g_{S(G^d)}-1$.
\end{prop}
\begin{proof}Let $S^\pm$ be surfaces obtained from $S(G^d)$ by a surgery
 on the curves of $\M^\pm$ respectively. Each $S^\pm$ is the inner boundary of the
 compression body $H^\pm$ together with several spheres that will be
 capped. Then the Euler characteristic of, say, $S^+$ is:
 \begin{eqnarray*}
 \chi(S^+)) = \chi((S(G^d)-\M^+) \cup 2g^+\ \mbox{discs}) =\\
 = \chi(S(G^d))+2g^+ = 2-2g_{S(G^d)}+2g^+
 \end{eqnarray*}
 This is also the sum of the Euler characteristics of the
 components $S_j$ of $S^+$,  each of which has genus $g_j$. Their number is exactly the
 number of the components of $S(G^d)-\M^+$. By Lemma
 \ref{separation} this number is $k^+$.
 Consequently,
 \[ 2-2g_{S(G^d)}+2g^+ = \sum_{j=1}^{k^+} (2-2g_j) = 2k^+-2\sum_{j=1}^{k^+} g_j=2k^+-2\pg^+,\]
 and the assertion follows.
\end{proof}
\begin{cor}\label{clmfldcond}
If $\pg^\pm=0$, then the manifold associated to $G^d$ is a closed
orientable ~3-manifold. If $\pg^+ = 0, \pg^- = 1$, or vice versa,
then the manifold associated to $G^d$ is a complement of a knot in
some closed orientable ~3-manifold.
\end{cor}
\begin{proof}
 Let $H^\pm$ be compression bodies obtained from $S(G^d) \times
 [-1,1]$ by gluing 2-handles along the curves of
 $\M^\pm$ respectively and filling in 3-balls.
 For a handlebody we need $\pg^\pm = 0$, since $\partial^-H^\pm$ has
 to be empty. To get a single torus as the inner
 boundary of $H^-$ and a single boundary component of $M$,
 we need $\pg^+ = 0,\pg^- = 1$.
\end{proof}
\begin{defn}
Two abstract Gauss diagrams will be called {\em equivalent} if one
can be obtained from the other by a sequence of $R,H,S,B$-moves
and their opposites.
\end{defn}
\begin{prop}
\label{pginvt} The numbers $\pg^\pm$ are invariants of the
equivalence classes of abstract Gauss diagrams.
\end{prop}
\begin{proof}
Consider a $B$-move. By Lemma \ref{bgs}, it preserves
$g_{S(G^d)}$. When a bubble is added to $\M^+$, a vertex is added
to both graphs $C^\pm(G^d)$. In $C^-$ this vertex is joined to
another by an edge, and $k^-$ does not change, but in $C^+$ no
edge joins it. Hence $k^+$ increases by $1$. The number of circles
also increases by $1$ only in $\M^+$. Computing $\pg^\pm$ by
formula of Proposition \ref{pg} shows that the values are
preserved by $B$-move.

Consider an $S$-move. Both $g^\pm$ increase by $1$. No new colors
are added, and the edges added are colored by $(c,c)$. Hence no
change is made on $k^\pm$. From Lemma \ref{sgs} we know that
$g_{S(G^d)}$ increases by $1$. Hence $\pg^\pm$ are preserved by
$S$-move.

Consider an $H$-move. This time $g_{S(G^d)}$ is preserved (Lemma
\ref{hgs}), and $g^\pm$ too. According to Lemma \ref{separation},
$k^\pm$ is the number of connected components of $S(G^d)-\M^\pm$.
Assume, for convenience, that the slide occurs in the $\M^+$
family. We would like to understand how it affects the connected
components of $S(G^d)-\M^\pm$. Since only the curves of $\M^+$ are
affected by the move, there is no change on the connected
components of $S(G^d)-\M^-$, and $k^-$ is preserved.

Let $\mu_1$ slide along $\mu_2$. Both of them are boundary
components of some $\Sigma'$, which is a connected component of
$S(G^d)-\M^+$. A handle slide looks on $\Sigma'$ like replacing
two components of $\partial\Sigma'$ by their connected sum. Assume
that $\Sigma'$ is glued along $\mu_1$ to another component
$\Sigma''$ of $S(G^d)-\M^+$. A slide of $\mu_1$ along $\mu_2$
looks like the addition to $\Sigma''$ of an annulus on a thin
band. No new components are created in $S(G^d)-\M^+$, and no
components are glued together. See Figure \ref{surfslide}.
\begin{figure}[ht]
\centerline{\psfig{figure=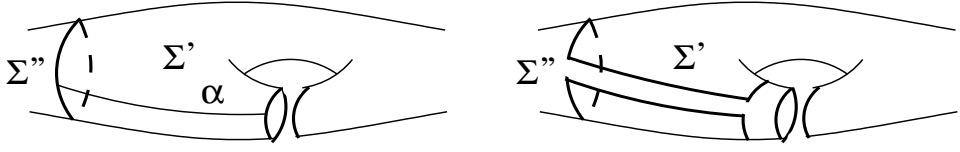,width=5in,silent=}}
\caption{Handle slide on the splitting surface} \label{surfslide}
\end{figure}

If both copies of $\mu_1$ are boundary components of $\Sigma'$,
then an annulus on a thin band is glued to the same $\Sigma'$.
There is no change in the other connected components of
$S(G^d)-\M^+$. Hence $k^+$ is preserved in both cases. We conclude
that $\pg^\pm$ is invariant under $H$-move.

Consider an $R$-move. There is no change in $g^\pm$. If the edges,
where the intersections appear, belonged to different cycles, then
the number of colors increases by one, and the number of cycles
does not change. This means that the genus of the surface
 associated to the changed diagram ${G'}^d$ is
$g_{S({G'}^d)}=1+\Delta_c-1-\frac{1}{2}(|h|+2-|c|)=g_{S(G^d)}$. In
this case the only new color, in both $C^\pm$, is a vertex of
valency $1$, joining the existing components. See also Figure
\ref{reidglue}. There is no other change in the graphs $C^\pm$. If
the edges where the intersections appeared belonged to the same
cycle, then the number of colors increases by $2$. The number of
cycles $|c|$ also increases by $2$. This means that
$g_{S({G'}^d)}=1+\Delta_c-\frac{1}{2}(|h|+2-|c|-2)=g_{S(G^d)}$.
One of the new colors, in both $C^\pm$, is a vertex of valency
$1$, joining the existing components, as in the previous case.
There are two more new colors replacing one old, and this looks
like splitting the vertex. Just as in the case of handle slide,
these vertices belong to the same component, since in both graphs
a path can be found between them. We conclude that $k^\pm$ is also
preserved and $\pg^\pm$ is invariant under $R$-move.
\end{proof}
\begin{them}
Closed connected orientable manifolds are in one-to-one
correspondence with the equivalence classes of abstract decorated
Gauss diagrams satisfying $\pg^\pm=0$.
\end{them}
\begin{proof}
 Obviously, two equivalent diagrams have the same associated
manifold. Now assume that we have a Gauss diagram $G^d_1$,
satisfying $\pg^\pm=0$, that is equivalent to a Gauss diagram
$G^d_2$ with $k^\pm=1$. By Proposition \ref{pginvt},
$g^\pm=g_{S(G^d_2)}$. By Corollary \ref{clmfldcond}, the manifold
associated to $G^d_2$ is closed. By Theorem \ref{gdthm} we are
done. It remains to show that any abstract Gauss diagram $G^d_1$
satisfying $\pg^\pm=0$ is equivalent to a Gauss diagram $G^d_2$
with $k^\pm=1$.

We can assume that $k^+> 1$. Let $\Sigma'$ be a connected
component of $S(G^d_1)-\M^+$. If $\Sigma'$ is not a disc with
holes, then $\Sigma'\times\{-1\}$ together with discs glued to its
boundary components is a non-trivial component of $\partial^-H^+$.
This contradicts the assumption $\pg^+=0$. Thus $\Sigma'$ is a
disc with holes.

Since $k^+>1$, there is a circle $\mu^+$ in $\partial\Sigma'$,
that is not glued to any other component of $\partial\Sigma'$.
Isotope $\mu^+$, probably performing some $R$-moves, so that
$\mu^+$ can be connected by arcs disjoint from $\M^-$ to all other
components of $\partial\Sigma'$. Then $\mu^+$ can slide along all
boundary components of $\Sigma'$. A handle slide looks on
$\Sigma'$ like replacing two components of $\partial\Sigma'$ by
their connected sum, hence the number of components of
$\partial\Sigma'$ decreases. From the proof of Proposition
\ref{pginvt} we know that handle slides do not change $k^\pm$.

The handle slides of $\mu^+$ along the components of
$\partial\Sigma'$ change $\Sigma'$ into a disc. Since
$\M^-\cap\Sigma'$ consists of disjoint arcs and bubbles, all these
can be removed by $B^{-1}$ and $R^{-1}$-moves. This will change
$\mu^+$ into a bubble. Remove the bubble by another $B^{-1}$-move.
This will decrease $k^+$ by $1$. Thus $G^d_1$ is equivalent to a
diagram with $k^+=1$.
\end{proof}

\subsection{The fundamental group and the first homology group}

 Since we reconstruct $M$
from its Gauss diagram as a cell complex, there is an obvious way
to obtain its fundamental group. A Gauss diagram is a graph, which
we denote also by $G$ for convenience.

\begin{prop} Let $G$ be an abstract Gauss diagram with all cycles colored
differently and $M$ the associated manifold. Let $T$ be the
maximal tree subgraph of $G$, containing all chords. The
generators of $\pi_1(M)$ are the edges of $G-T$, and the relations
are defined by the cycles of $G$ and the circles of both families.
\end{prop}
\begin{proof}
The fundamental group of the ribbon graph $\Gamma(G)$ coincides
with the fundamental group of $G$, which is the free group
generated by the edges of $G-T$. To construct the splitting
surface one glues discs to the boundary components of $\Gamma(G)$,
which means to the cycles. For any cycle write a relation going
along it and writing down only the edges in $G-T$. In the cycles
part of the edges are traversed in the direction opposite to their
orientation. In this case write the inverse of the generator.

Gluing 2-handles along the circles of $G$ means introducing
another set of relations. The circles are also sequences of edges,
and only those in $G-T$ should be written. Filling in 3-balls does
not affect the fundamental group.
\end{proof}
\begin{note}
There is another, simpler way to obtain the fundamental group for
closed orientable manifolds, if the number of circles in each
family of the Gauss diagram equals the genus of the associated
surface and no separation occurs. This is the method of
Proposition \ref{pi1}.
\end{note}

 The first homology group of $M$ is the abelian group generated by
 the edges of $G-T$. The relations coming from the circles in the Gauss
diagram remain unchanged. In some cycles there could be edges
appearing twice, with opposite signs. Hence in $H_1(M)$ in the
relations coming from the cycles one writes only those generators
that are colored by two different cycles.

\subsection{Examples}
In both following examples we represent the diagrams with all
cycles colored differently, hence the decorations are omitted.
\begin{description}
\item[Thickened torus] $T^2\times I$ can be represented by a
diagram in Figure \ref{torusxI}, which also shows how the curves
are situated on the surface. The count of cycles gives $|c|=8$, so
$g_{S(G)}=1$, which equals the number of curves in each family.
One can also check, using the graphs $C^\pm$ of Section
\ref{bdry}, that each curve separates the surface into two
components. The Proposition \ref{pg} then implies that both
boundary components are tori.
\begin{figure}[ht]
\centerline{\psfig{figure=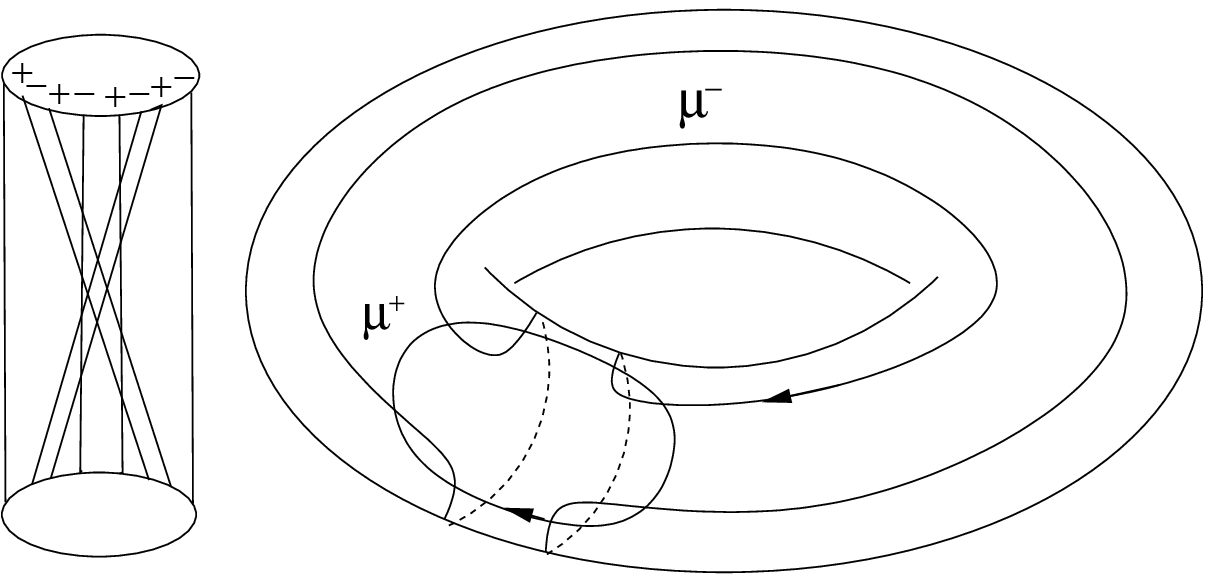,height=1.2in,silent=}}
\caption{$T^2\times I$} \label{torusxI}
\end{figure}

\item[Unknot complement] The diagram of Figure \ref{unknot} has
$4$ cycles, which implies genus $g=1=g_{S(G)}$, and drawing the
graphs $C^\pm$ one can see that the top family does not separate
the surface, while the bottom does, and we have two components.
Let us count the fundamental group. The maximal tree $T$ is
thickened in the diagram. The generators are the arcs in $G-T$, as
shown in Figure \ref{unknot}. The fundamental group then is:
\[
\langle e_1,e_2,e_3,g_1,g_2|e_1e_2e_3,g_1g_2,
     e_1g_1^{-1}e_1^{-1}g_2^{-1},e_2g_1,e_2^{-1}e_3g_2,e_3^{-1}\rangle,
\]
where the last four relations come from the cycles. This gives us
$\pi_1(M)=Z$. The manifold is the complement of the unknot in
$S^3$.
\begin{figure}[ht]
\centerline{\psfig{figure=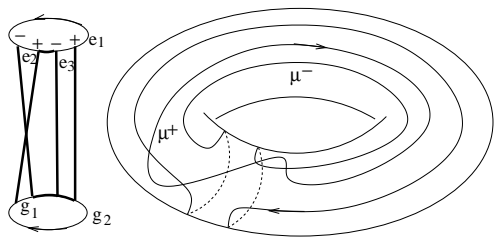,height=1.2in,silent=}}
\caption{The complement of unknot in $S^3$} \label{unknot}
\end{figure}
\end{description}

\end{document}